\documentclass[11pt]{amsart}
\usepackage{fullpage}
\usepackage{amssymb}
\usepackage{mathtools,color}
\usepackage{amsrefs}
\usepackage{leftidx}
\usepackage{subcaption}
\usepackage{tikz}
\usepackage{eucal}
\usepackage{appendix}
\usepackage[T1]{fontenc}

\newcommand{\T}{\mathcal{T}}
\renewcommand{\emptyset}{\text{\O}}
\newcommand{\CC}{\mathcal{C}}
\newcommand{\N}{\mathbb{N}}
\newcommand{\la}{\lambda}
\renewcommand{\SS}{\mathbf{S}}

\theoremstyle{plain}
\newtheorem{theorem}{Theorem}[section]
\newtheorem{prop}[theorem]{Proposition}
\newtheorem{lemma}[theorem]{Lemma}

\theoremstyle{definition}
\newtheorem{rem}[theorem]{Remark}
\newtheorem{example}[theorem]{Example}

\numberwithin{equation}{section}

\allowdisplaybreaks

\begin{document}

\title{The sum of all width-one matrices}

\author{William Q. Erickson}
\address{
William Q.~Erickson\\
Department of Mathematics\\
Baylor University \\ 
One Bear Place \#97328\\
Waco, TX 76798} 
\email{Will\_Erickson@baylor.edu}

\author{Jan Kretschmann}
\address{
Jan Kretschmann\\
Department of Mathematical Sciences\\
University of Wisconsin--Milwaukee \\ 
3200 N.~Cramer St.\\
Milwaukee, WI 53211} 
\email{kretsc23@uwm.edu}

\begin{abstract}
A nonnegative integer matrix is said to be width-one if its nonzero entries lie along a path consisting of steps to the south and to the east. 
These matrices are important in optimal transport theory: the northwest corner algorithm, for example, takes supply and demand vectors and outputs a width-one matrix.
The problem in this paper is to write down an explicit formula for the sum of all width-one matrices (with given dimensions $n \times n$ and given sum $d$ of the entries).
We prove two strikingly different formulas.  
The first, a $\leftidx{_4}{F}{_3}$ hypergeometric series with unit argument, is obtained by applying the Robinson--Schensted--Knuth correspondence to the width-one matrices; the second is obtained via Stanley--Reisner theory.
Computationally, our two formulas are complementary to each other: the first formula outperforms the second if $d$ is fixed and $n$ increases, while the second outperforms the first if $n$ is fixed and $d$ increases.
We also show how our result yields a new non-recursive formula for the mean value of the discrete earth mover's distance (i.e., the solution to the transportation problem), whenever the cost matrix has the Monge property.
\end{abstract}



\maketitle

\section{Introduction}

Let $\T(d,n)$ be the set of all $n \times n$ matrices with nonnegative integer entries summing to $d$, such that the nonzero entries lie along a path consisting of steps to the south and to the east.
For example, one element of $\T(30,5)$ is the matrix
\[
\left[\begin{smallmatrix}
    5 & 3 & \textcolor{gray!50}{0} & \textcolor{gray!50}{0} & \textcolor{gray!50}{0} \\
    \textcolor{gray!50}{0} & 2 & \textcolor{gray!50}{0} & \textcolor{gray!50}{0}  & \textcolor{gray!50}{0} \\
    \textcolor{gray!50}{0} & \textcolor{gray!50}{0}  & 3 & \textcolor{gray!50}{0}  & \textcolor{gray!50}{0} \\
    \textcolor{gray!50}{0} & \textcolor{gray!50}{0}  & \textcolor{gray!50}{0}  & \textcolor{gray!50}{0}  & \textcolor{gray!50}{0} \\
    \textcolor{gray!50}{0} & \textcolor{gray!50}{0} & 9 & 7 & 1
\end{smallmatrix}
\right].
\]
We call these \emph{width-one} matrices, since their support has width 1 in the underlying poset of matrix coordinates; 
moreover, they are the exponent matrices of width-one monomials, in the language of Sturmfels~\cite{Sturmfels}*{p.~138}.
In this paper, we solve the problem of writing down an explicit formula for the sum $\SS(d,n)$ of all matrices in $\T(d,n)$.

\begin{figure}[t]

    \centering
    \begin{subfigure}[b]{0.45\textwidth}
         \centering
         \includegraphics[height=.7\textwidth,width=\textwidth]{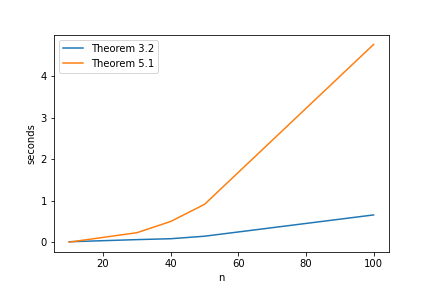}
         \caption{}
         \label{subfig:CompA}
    \end{subfigure}
    \begin{subfigure}[b]{0.45\textwidth}
        \centering
        \includegraphics[height=0.7\textwidth,width=\textwidth]{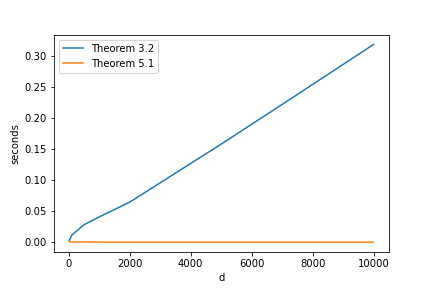}
        \caption{}
        \label{subfig:CompB}
    \end{subfigure}

    \caption{Comparison of computing time with respect to the parameters $d$ and $n$. 
    In \ref{subfig:CompA}, we fix $d=30$ and compare the runtime (in seconds) of both approaches for varying $n$. 
    In \ref{subfig:CompB}, we fix $n=5$ and let $d$ vary. 
    Some values of $\SS(d,5)$ can be seen in Appendix \ref{sec:appValues}.}
    \label{fig:complexity}
\end{figure}

In fact, we present two different formulas for the $(i,j)$ entry of $\SS(d,n)$.
On one hand (Theorem~\ref{thm:RSK}), we show that this entry equals $\sum_{k=0}^{d-1} \binom{i + k-1}{k} \binom{j + k-1}{k}\binom{n - i + d -k - 1}{n-i}  \binom{n - j + d-k-1}{n-j}$, more succinctly expressed as
\[
\binom{i+d-2}{d-1}\binom{j+d-2}{d-1} \leftidx{_4}{F}{_3}\!\left[\genfrac{}{}{0pt}{0}{n-i+1,\: n-j+1, \: 1-d, \: 1-d}{1, \: 2-d-i, \: 2-d-j} \: ; \: 1 \right],
\]
where $\leftidx{_4}{F}{_3}$ is the generalized hypergeometric series evaluated at unity.
On the other hand (Theorem~\ref{thm:main result h-polys}), this $(i,j)$ entry also equals the sum
\[
    \sum_{k=0}^{\min\{d,n\}-1} \binom{2n+d-k-2}{2n-1} \sum_{\ell = 0}^k \binom{i-1}{\ell}\binom{j-1}{\ell}\binom{n-i}{k-\ell}\binom{n-j}{k-\ell}.
\]
The equivalence of these formulas is interesting in its own right, but we also include them both because of the computational difference between them; see Figure~\ref{fig:complexity}.

As for methodology, the first formula actually follows from a straightforward application of the famous Robinson--Schensted--Knuth (RSK) correspondence; we give the proof in Section~\ref{sec:RSK}.
Our second formula, on the other hand, is obtained by translating the problem into combinatorial commutative algebra.
Hence we first devote Section~\ref{sec:ASC} to the necessary theory of simplicial complexes and Stanley decompositions.
Our primary tool is a certain shelling of the order complex on the poset of matrix coordinates, described in detail in the paper~\cite{BCS} by Billera--Cushman--Sanders on the harmonic oscillator.
This shelling induces a Stanley decomposition of the Stanley--Reisner ring of the order complex, and also easily yields the $h$-polynomial of certain subcomplexes.
We then obtain our second formula in Section~\ref{sec:main result h-poly}, by isolating the exponent of a single variable $x_{ij}$ in the product of all monomials inside the $d$th graded component of the Stanley--Reisner ring.

Finally, in Section~\ref{sec:EMD} we point out one application of our result to optimal transport theory (which did in fact motivate the problem in this paper).
There is a well-known algorithm in the transport theory literature known as the ``northwest corner rule,'' which associates a pair of supply and demand vectors to a unique width-one matrix $T \in \T(d,n)$. 
(This matrix is called a \emph{transport plan}, which motivates the notation $\T$ and $T$ in the paper.)
Whenever the cost of transport satisfies the Monge property~\eqref{Monge property}, this $T$ gives the solution (often called the \emph{earth mover's distance}, or EMD) to the classical transportation problem: in other words, $T$ encodes directions for transporting material from the supply locations to the demand locations such that the total cost is minimized.
It is natural to seek the mean value of the EMD, taken over all possible supply and demand vectors; we show that our formulas for $\SS(d,n)$ immediately yield this mean value.
As far as we are aware, previous methods for computing the mean EMD have been recursive in nature, and limited to particular cost functions; by contrast, the combinatorial result in this paper offers a direct and far more flexible method for computing this mean value.

\section{Notation and hypergeometric series}

\subsection{Detailed statement of the problem}

For positive integers $i$ and $j$, we define the following poset, equipped with the product order $\preceq$:
\begin{equation}
    \label{def:Pi}
    \Pi_{ij} \coloneqq \{(x,y) : 1 \leq x \leq i, \: 1 \leq y \leq j\}, \qquad
    (x,y) \preceq (x',y') \iff x \leq x' \text{ and } y \leq y'.
\end{equation}
Elsewhere in the literature, $\Pi_{ij}$ is sometimes written as $[i] \times [j]$.
We will abbreviate $\Pi_n \coloneqq \Pi_{nn}$.  
Because our main problem deals with $n \times n$ matrices, we will always be working inside $\Pi_{n}$, but it will be useful to consider its subposets $\Pi_{ij}$, which are the lower-order ideals generated by each matrix position $(i,j)$.

We adopt the standard terminology from order theory: in particular, a \emph{chain} is a totally ordered subset of $\Pi_{ij}$, and an \emph{antichain} is a subset whose elements are pairwise incomparable. 
The \emph{width} of a subset $S \subseteq \Pi_{ij}$ is the size of the largest antichain contained in $S$. 
Equivalently, by Dilworth's theorem, the width of $S$ equals the minimum number of chains into which $S$ can be partitioned.  

The \emph{support} of an $n \times n$ matrix $M$ is the set
\[
\operatorname{supp}(M) \coloneqq \{(i,j) \in \Pi_{n} : M_{ij} \neq 0\}.
\] 
As mentioned above, we say that $M$ is a \emph{width-one} matrix if $\operatorname{supp}(M) \subseteq \Pi_{n}$ has width 1.  
Note that when $\Pi_{n}$ is visualized as the coordinates of an $n \times n$ matrix, the relation $\preceq$ means ``weakly northwest''; therefore, in a width-one matrix, all nonzero entries lie along a path from the northwest to the southeast corner (i.e., a maximal chain), consisting of steps to the south and to the east. 
Now let ${\rm M}_n(\N)$ denote the set of $n \times n$ matrices with entries in $\N \coloneqq \{0, 1, 2, \ldots\}$. For a positive integer $d$, we define the set
\begin{equation}
    \label{T(d,n) definition}
    \T(d,n) \coloneqq \{ T \in {\rm M}_n(\N) : \text{$T$ is width-one and $\sum_{i,j} T_{ij} = d$}\}.
\end{equation} 
The main problem of this paper is to write down an explicit formula for the sum of all these matrices, which we call
\begin{equation}
    \label{Bold T definition}
    \SS(d,n) \coloneqq \sum_{\mathclap{T \in \T(d,n)}} T.
\end{equation}

\subsection{Hypergeometric series}

First, it will be a useful fact that the number of weakly increasing sequences of length $a$, with entries from $\{1, \ldots, k\}$, equals the number of weak integer compositions of $a$ into $k$ parts, which is well known to be the binomial coefficient
\begin{equation}
\label{composition formula}
    \binom{a + k - 1}{a} = \binom{a+k - 1}{k-1} \coloneqq \frac{(a+k-1)!}{a!(k-1)!}.
\end{equation}
Next, we recall the (ordinary, or Gaussian) \emph{hypergeometric series}
\[
    \leftidx{_2}{F}{_1}\!\left[\genfrac{}{}{0pt}{0}{a,b}{c} ; z\right] \coloneqq \sum_{k=0}^\infty \frac{(a)_k (b)_k}{(c)_k} \frac{z^k}{k!},
\]
where $(a)_k \coloneqq a(a+1)\cdots(a+k-1)$ is the Pochhammer symbol for the rising factorial.
In the special case $c=1$, it is easy to see that
\begin{equation}
    \label{hypergeo special case}
    \leftidx{_2}{F}{_1}\!\left[\genfrac{}{}{0pt}{0}{a,b}{1} ; z\right] = \sum_{k=0}^\infty \frac{(a)_k}{k!} \frac{(b)_k}{k!} z^k = \sum_{k=0}^\infty \binom{a+k-1}{k} \binom{b+k-1}{k} z^k,
\end{equation}
where the second equality holds only if $a,b \in \N$.
The \emph{generalized} hypergeometric series $\leftidx{_p}{F}{_q}$ is defined analogously:
\[
\leftidx{_p}{F}{_q}\!\left[\genfrac{}{}{0pt}{0}{a_1, \ldots, a_p}{b_1, \ldots, b_q}; z\right] \coloneqq \sum_{k=0}^\infty \frac{(a_1)_k \cdots (a_p)_k}{(b_1)_k \cdots (b_q)_k} \frac{z^k}{k!}.
\]
Since our interest in this paper is purely combinatorial, we disregard issues of convergence and treat $\leftidx{_p}{F}{_q}$ as a formal power series.
Note that if some $b_j$ is a nonpositive integer, then infinitely many of the coefficients are undefined.
But if, for example, some $a_i$ is also a nonpositive integer with $a_i \geq b_j$, then the series terminates before these undefined coefficients, and so the series as a whole is still defined (and is a polynomial in $z$).

The discovery of identities involving hypergeometric series is a longstanding, and yet still quite active, research area within combinatorics; see the books \cite{Erdelyi}*{Ch.~II} and \cite{PWZ}, for example.
We will appeal to the following identity 4.3(14) in~\cite{Erdelyi}*{p.~187}, which expresses the convolution of two ordinary hypergeometric series:
\[
\leftidx{_2}{F}{_1}\!\left[\genfrac{}{}{0pt}{0}{a,b}{c}; \alpha z\right] \leftidx{_2}{F}{_1}\!\left[\genfrac{}{}{0pt}{0}{a',b'}{c'}; \beta z\right] = \sum_{k=0}^\infty \frac{(a)_k (b)_k}{(c)_k}\frac{(\alpha z)^k}{k!} \leftidx{_4}{F}{_3}\!\left[\genfrac{}{}{0pt}{0}{a', \:b', \: 1-k-c, \: -k}{c', \: 1-k-a, \: 1-k-b}; \beta/\alpha\right].
\]
We will encounter the specialization where $c=c'=\alpha=\beta=1$, which yields
\begin{equation}
\label{sum of 2F1s}
\leftidx{_2}{F}{_1}\!\left[\genfrac{}{}{0pt}{0}{a,b}{1}; z\right] \leftidx{_2}{F}{_1}\!\left[\genfrac{}{}{0pt}{0}{a',b'}{1}; z\right] = \sum_{k=0}^\infty 
\underbrace{\frac{(a)_k}{k!} \frac{(b)_k}{k!}}_{\substack{\binom{a+k-1}{k} \binom{b+k-1}{k} \\[.75ex] \text{if $a,b \in \N$}}} \leftidx{_4}{F}{_3}\!\left[\genfrac{}{}{0pt}{0}{a', \: b', \: -k, \:-k}{1, \:1-k-a, \:1-k-b}; 1\right] z^k.
\end{equation}

\section{Main result, first version}
\label{sec:RSK}

Any width-one matrix $T \in \T(d,n)$ can be written as a sum of elementary matrices
\[
T = \sum_{i,j=1}^n T_{ij} E_{ij} = \sum_{\mathclap{(i,j) \in {\rm supp}(T)}} T_{ij} E_{ij},\]
where $E_{ij}$ is the $n \times n$ matrix whose $(i,j)$ entry is $1$, with $0$'s everywhere else.
Then by definition, the coordinates $(i,j)$ appearing in the right-hand sum form a chain in $\Pi_n$.
Writing out this chain of $(i,j)$'s as columns $\begin{smallmatrix}i\\j\end{smallmatrix}$ in ascending order, where each $\begin{smallmatrix}i\\j\end{smallmatrix}$ occurs $T_{ij}$ times, we obtain a $2 \times d$ array, often called a \emph{biword} in the literature:
\begin{equation}
    \label{biword}
\begin{pmatrix}
    i_1 & i_2 & \ldots & i_d\\
    j_1 & j_2 & \ldots & j_d
\end{pmatrix}, \quad \text{each row weakly increasing with entries in $\{1, \ldots, n\}$}.
\end{equation}
For example, taking a matrix in $\T(10,4)$, we obtain the following biword:
\[
\left[\begin{smallmatrix}
    \textcolor{gray!50}{0} & 3 & \textcolor{gray!50}{0} & \textcolor{gray!50}{0}\\
    \textcolor{gray!50}{0} & 2 & \textcolor{gray!50}{0} & \textcolor{gray!50}{0}\\
    \textcolor{gray!50}{0} & 1 & 2 & \textcolor{gray!50}{0}\\
    \textcolor{gray!50}{0} & \textcolor{gray!50}{0} & 1 & 1
    \end{smallmatrix}
    \right] 
    \leadsto
    \begin{pmatrix}
    1&1&1&2&2&3&3&3&4&4\\
    2&2&2&2&2&2&3&3&3&4
    \end{pmatrix}.
\]
This procedure is clearly invertible: given any biword of the form~\eqref{biword}, we recover the corresponding width-one matrix $T \in \T(d,n)$ by taking the sum $\sum_{\ell = 1}^d E_{i_\ell, j_\ell}$.
The upshot is that we have a bijection between $\T(d,n)$ and the set of biwords of the form~\eqref{biword}.

\begin{rem}
The reader may recognize this bijection as a very special case of the Robinson--Schensted--Knuth (RSK) correspondence, restricted to width-one matrices.
In this case, the two corresponding semistandard Young tableaux have only one row each, and are precisely the rows of the biword~\eqref{biword}.
This follows directly from Knuth's construction in~\cite{Knuth}*{\S3}; in general, the width of the support of a matrix equals the number of rows in the corresponding tableaux.
\end{rem}

\begin{theorem}
    \label{thm:RSK}
    For all $1 \leq i,j \leq n$, the $(i,j)$ entry of $\SS(d,n)$ is
    \[
    \SS(d,n)_{ij} = \binom{i+d-2}{d-1}\binom{j+d-2}{d-1} \leftidx{_4}{F}{_3}\!\left[\genfrac{}{}{0pt}{0}{n-i+1,\: n-j+1, \: 1-d, \: 1-d}{1, \: 2-d-i, \: 2-d-j} \: ; \: 1 \right].    
    \]
\end{theorem}

\begin{proof}
    For each $T \in \T(d,n)$, consider its corresponding biword~\eqref{biword}, in which each column $\begin{smallmatrix}i\\j\end{smallmatrix}$ contributes 1 to the $(i,j)$ entry in $T$.
    Hence the entry $\SS(d,n)_{ij}$ equals the total number of occurrences of the column $\begin{smallmatrix}i\\j\end{smallmatrix}$ within all possible biwords.

    Suppose that the $\ell$th column of a biword is $\begin{smallmatrix}i\\j\end{smallmatrix}$.
    In the top row, this implies that the $\ell-1$ entries to the left of $i$ lie in the set $\{1,\ldots, i\}$, and the $d-\ell$ entries to the right of $i$ lie in the set $\{i, i+1, \ldots, n\}$, which contains $n + 1 - i$ elements.
    Hence by~\eqref{composition formula}, the following product of binomial coefficients equals the number of ways to fill the top row of the biword such that the $\ell$th entry is $i$:
    \[
    \binom{(\ell-1) + i - 1}{\ell-1} \binom{(d-\ell) + (n + 1 - i) - 1}{d-\ell}.
    \]
    Using the same argument for the bottom row of the biword (replacing $i$ by $j$), and then multiplying the top and bottom results, we conclude that the number of biwords with $\begin{smallmatrix}i\\j\end{smallmatrix}$ as the $\ell$th column equals
    \[
    \binom{i + \ell - 2}{\ell-1} \binom{j + \ell - 2}{\ell-1} \binom{n - i + d - \ell}{d - \ell}  \binom{n - j + d - \ell}{d - \ell}.
    \]
    To obtain the number of times $\begin{smallmatrix}i\\j\end{smallmatrix}$ occurs as any column in a biword, we sum over all columns $\ell = 1, \ldots, d$.
    Following this with the substitution $k = \ell-1$, we have
   \[
    \SS(d,n)_{ij} = \sum_{k=0}^{d-1} \underbrace{\binom{i + k-1}{k} \binom{j + k-1}{k}}_{
    \substack{
    \text{coeff. of $z^k$ in}
    \\[.5ex]
    \leftidx{_2}{F}{_1}\left[\genfrac{}{}{0pt}{2}{i,j}{1} ; z\right] 
    }
    } \underbrace{\binom{n - i + d - 1 -k}{d - 1 - k}  \binom{n - j + d-1-k}{d-1-k}}_{\substack{
    \text{coeff. of $z^{(d-1)-k}$ in}
    \\[.5ex]
    \leftidx{_2}{F}{_1}\left[\genfrac{}{}{0pt}{2}{n-i+1,\:n-j+1}{1} ; z\right]
    }},
   \]
where we have recognized the two coefficients from the hypergeometric series in~\eqref{hypergeo special case}.
This makes it clear that $\SS(d,n)_{ij}$ is the coefficient of $z^{d-1}$ in the product
\begin{align*}
    &\leftidx{_2}{F}{_1}\!\left[\genfrac{}{}{0pt}{0}{i,j}{1} ; z\right]\: \leftidx{_2}{F}{_1}\!\left[\genfrac{}{}{0pt}{0}{n-i+1,\: n-j+1}{1} ; z\right]\\[2ex]
    =& \sum_{k=0}^\infty 
\binom{i+k-1}{k} \binom{j+k-1}{k} \leftidx{_4}{F}{_3}\!\left[\genfrac{}{}{0pt}{0}{n-i+1, \: n-j+1, \: -k, \:-k}{1, \:1-k-i, \:1-k-j}; 1\right] z^k,
\end{align*}
where the equality is just the identity~\eqref{sum of 2F1s}.
Reading off the coefficient for $k = d-1$, we obtain the expression in the theorem.
\end{proof}

\begin{figure}[ht]

    \centering
    \begin{subfigure}[b]{0.45\textwidth}
         \centering
         \includegraphics[height=.8\textwidth,width=.8\textwidth]{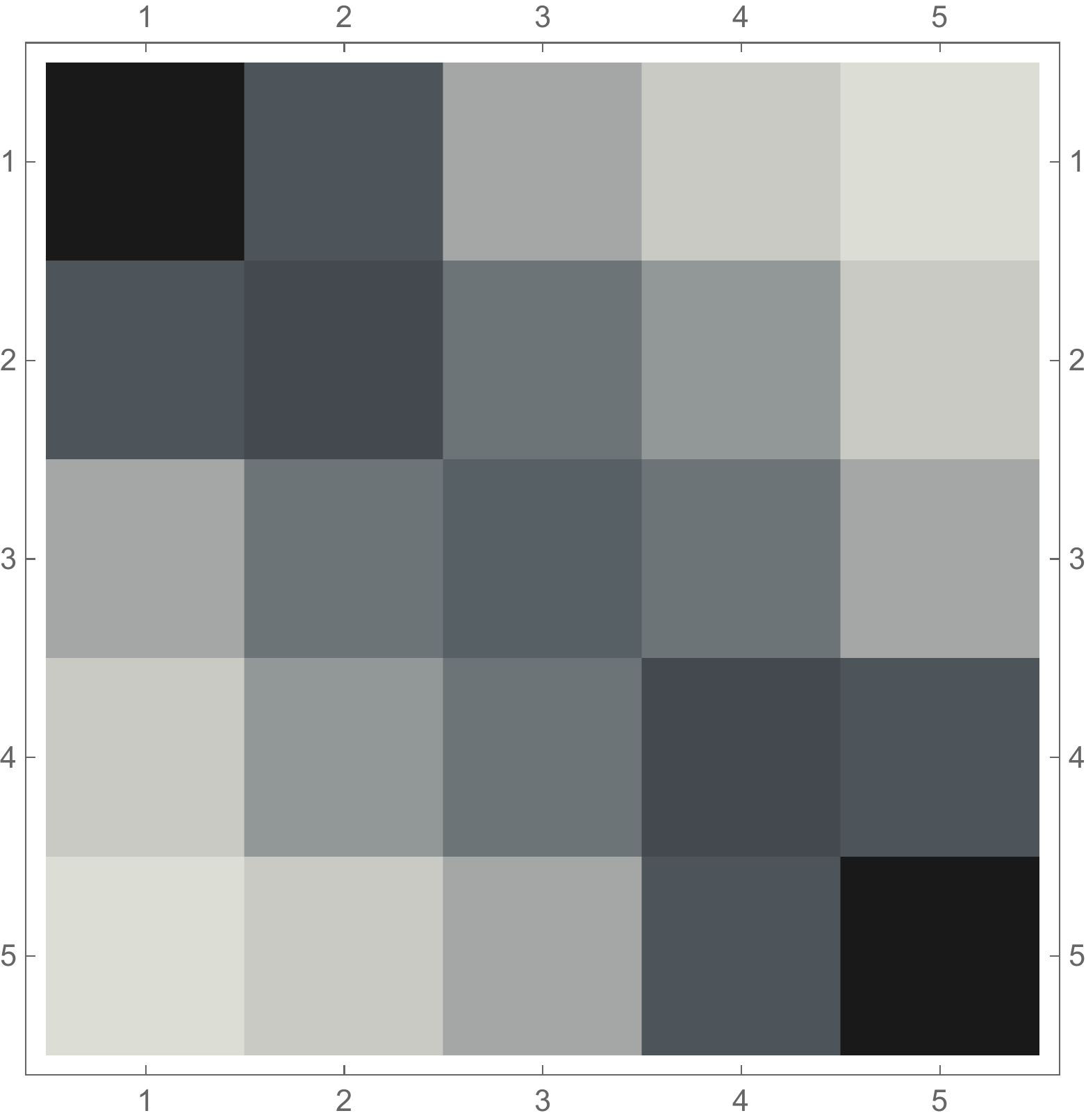}
         \caption{Matrix plot of $\SS(30,5)$.}
    \end{subfigure}
    \begin{subfigure}[b]{0.45\textwidth}
        \centering
        \includegraphics[height=0.8\textwidth,width=0.8\textwidth]{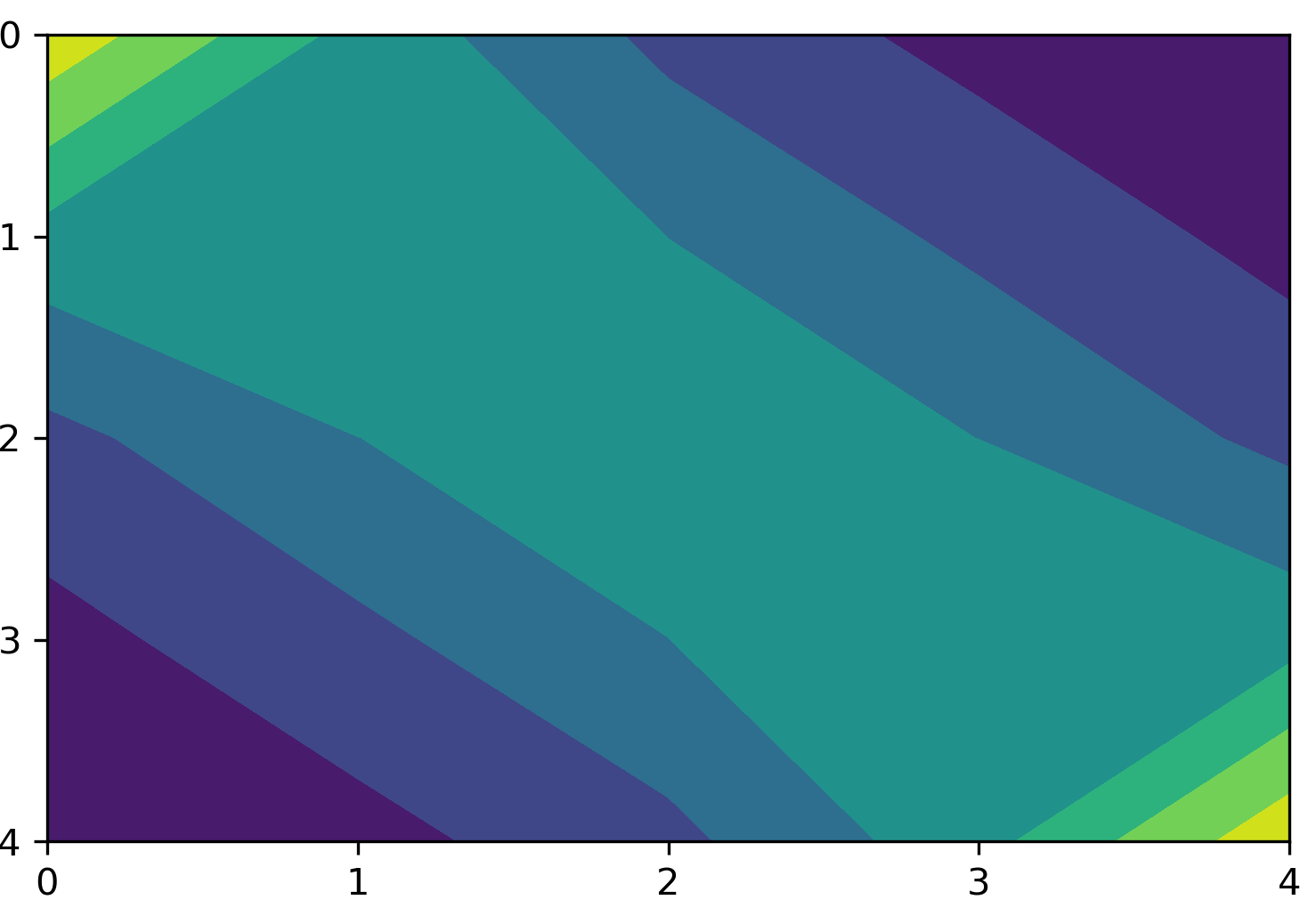}
        \caption{Contour plot of $\SS(30,5)$.}
    \end{subfigure}

\vspace{0.5cm}

    \begin{subfigure}[b]{0.45\textwidth}
         \centering
         \includegraphics[height=.8\textwidth,width=.8\textwidth]{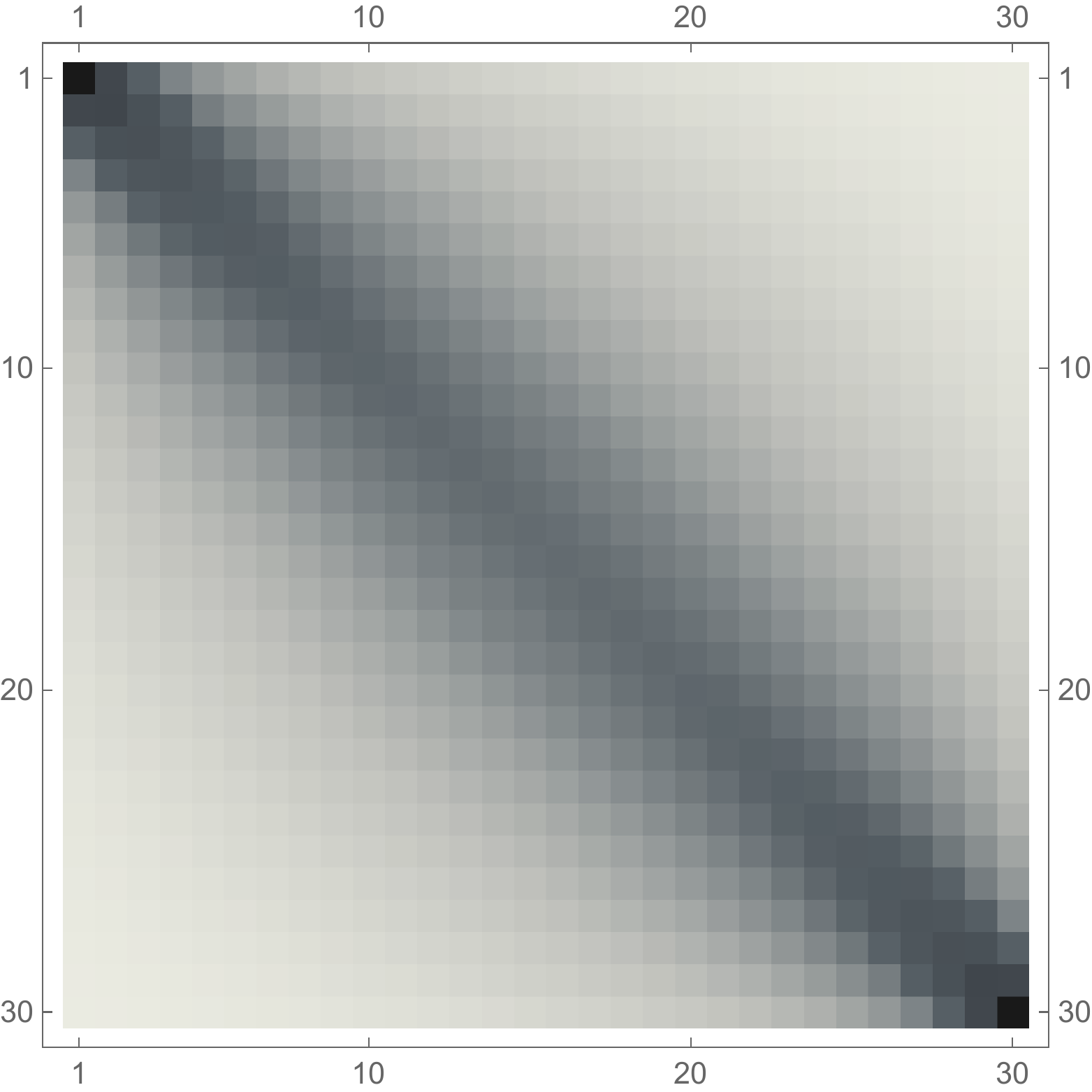}
         \caption{Matrix plot of $\SS(10000, 30)$.}
    \end{subfigure}
     \begin{subfigure}[b]{0.45\textwidth}
         \centering
         \includegraphics[height=.8\textwidth,width=.8\textwidth]{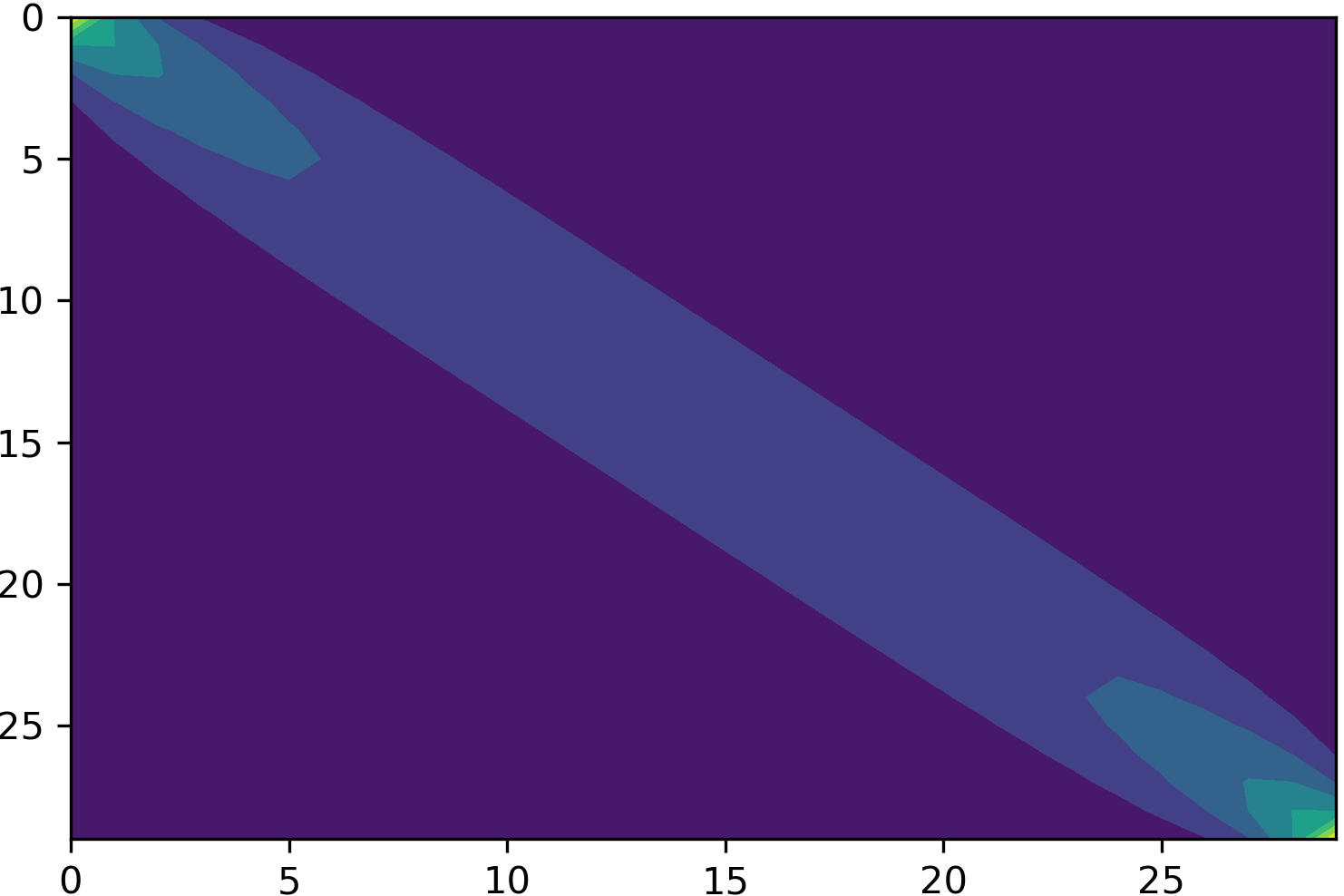}
         \caption{Contour plot of $\SS(10000, 30)$.}
     \end{subfigure}

    \caption{Matrix and contour plots
    of $\SS(d,n)$.
    In the matrix plots, darker shades correspond to greater entries; in the contour plots, blue and yellow correspond to the lowest and highest level curves, respectively.
    In $\SS(30,5)$, the entries lie in the approximate interval $[1.6 \times 10^8, \: 6.7 \times 10^9]$.
    In $\SS(10000, 30)$, the entries lie in the approximate interval $[8.6 \times 10^{155}, \: 2.4 \times 10^{172}]$.}
    \label{fig:matrix plot}
\end{figure}

As an example, Figure~\ref{fig:matrix plot} shows the matrix plot and contour plot of $\SS(30,5)$ and $\SS(10000,30)$.  
Due to the obvious symmetries in $\T(d,n)$ given by reflecting width-one matrices about the main diagonal and the anti-diagonal, it is clear that $\SS(d,n)$ is bisymmetric;
therefore in practice, we need only compute $\SS(d,n)_{ij}$ for $i \leq j \leq n+1-i$.
In the remainder of the paper, we develop an argument using a Stanley--Reisner ring to obtain a second formula which even more efficiently handles cases where $d$ is much larger than $n$.

\begin{rem}
    The only reason we restrict our attention to \emph{square} width-one matrices is to keep notation uncluttered; it is trivial to adapt Theorem~\ref{thm:RSK} to matrices of arbitrary dimensions $n_1 \times n_2$:
    \[
        \SS(d, n_1 \times n_2)_{ij} = \binom{i+d-2}{d-1}\binom{j+d-2}{d-1}\leftidx{_4}{F}{_3}\!\left[\genfrac{}{}{0pt}{0}{n_1-i+1,\: n_2-j+1, \: 1-d, \: 1-d}{1, \: 2-d-i, \: 2-d-j} \: ; \: 1 \right].
    \]
    \label{rem:n1n2}
\end{rem}

\section{Simplicial complexes}
\label{sec:ASC}

This section marks our turning point toward combinatorial commutative algebra, in order to derive a second formula for $\SS(d,n)$. 
We largely follow the exposition in~\cite{StanleyAC}*{Ch.~12}; see also~\cite{24Hours}*{Ch.~16}.  
We outline the general theory in the first two subsections, before adapting it to our main problem in the third.
 
 \subsection{Abstract simplicial complexes}
 \label{sub:simplicial complexes}

 Given a finite set $V$, an \emph{(abstract) simplicial complex} on $V$ is a collection $\Delta$ of subsets of $V$ such that
 \begin{itemize}
     \item $\{v\} \in \Delta$ for all $v \in V$;
     \item if $S \in \Delta$ and $R \subseteq S$, then $R \in \Delta$.
 \end{itemize}
 Elements of $V$ are called \emph{vertices}, and $V$ is called the \emph{vertex set}.  The elements of $\Delta$ are called \emph{faces}, and the \emph{dimension} of a face is one less than its cardinality.  The dimension of $\Delta$ is the maximum of the dimensions of its faces.  
 
 Let $\Delta$ be a nonempty simplicial complex of dimension $m-1$. 
We denote by $f_i$ the number of faces of dimension $i$.  In particular, $f_{-1} = 1$ since $\emptyset \in \Delta$.  The \emph{$f$-vector} is the sequence $(f_0, \ldots, f_{m-1})$.
It will be more convenient for us to work instead with the \emph{$h$-vector} $(h_0, \ldots, h_m)$, where the numbers $h_\ell$ are defined as follows:
\begin{equation}
\label{eq-h-to-f}
h_\Delta(t) \coloneqq \sum_{\ell=0}^m h_\ell t^\ell = \sum_{i=0}^m f_{i-1} t^i (1-t)^{m-i}.
\end{equation}
We call $h_\Delta(t)$ the \emph{$h$-polynomial} of $\Delta$.
 
A maximal face in $\Delta$ (with respect to inclusion) is called a \emph{facet}. 
We say that $\Delta$ is \emph{pure} if every facet has the same dimension.  
A pure simplicial complex is said to be \emph{shellable} if there exists an ordering $F_1, \ldots, F_{s}$ of its facets with the following property: for all $i = 1, \ldots, s$, the power set of $F_i$ has a unique minimal element ${\rm R}(F_i)$ not belonging to the subcomplex generated by $F_1, \ldots, F_{i-1}$.  
Such an ordering is called a \emph{shelling}, and ${\rm R}(F_i)$ is called the \emph{restriction} of the facet $F_i$.  
Crucial to our method is the following combinatorial description of the $h$-polynomial: if $F_1, \ldots, F_s$ is a shelling of $\Delta$, then we have
\begin{equation}
    \label{h-vec as restrictions}
    h_\Delta(t) = \sum_{i=1}^s t^{\#{\rm R}(F_i)}.
\end{equation}
In other words, $h_\ell$ counts the number of facets whose restrictions have size $\ell$ (and this is independent of the choice of shelling):
\[ 
    h_\ell = \#\{i : \#{\rm R}(F_i) = \ell\}.
\]

\begin{example}
    Consider the simplicial complex that is the boundary of the octahedron (see Figure~\ref{fig:octahedron}).
    \begin{figure}[h]
        \centering
        \begin{tikzpicture}[z={(-.3cm,-.2cm)}, 
line join=round, line cap=round,
every node/.style = {inner sep=1pt, font=\scriptsize}
]
\draw ( 0,2,0) node [above] {$A$} --  (-2,0,0) node [left] {$B$} --  (0,-2,0) node [below] {$C$} --  ( 2,0,0) node [right] {$D$} --  ( 0,2,0) --  ( 0,0,2) node [below left] {$E$} --  (0,-2,0) (2,0,0) -- (0,0,2) -- (-2,0,0);
\draw [dashed] (0,2,0) -- (0,0,-2) node [above] {\;\;$F$} -- (0,-2,0) (2,0,0) -- (0,0,-2) -- (-2,0,0);
\end{tikzpicture}
        \caption{The octahedron.}
        \label{fig:octahedron}
    \end{figure}
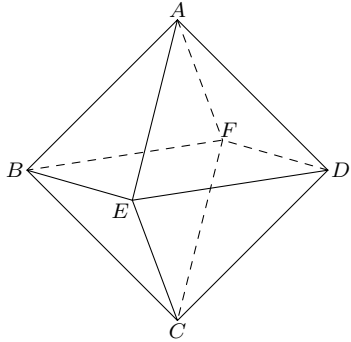

To show that this complex is shellable, it suffices to find an ordering $F_1,\ldots,F_8$ of its eight facets, such that each $F_i$ contains a unique minimal element ${\rm R}(F_i)$ that is not contained in the subcomplex generated by $F_1, \ldots, F_{i-1}$.
Note that when $i=1$, this subcomplex is the empty set.
Below, we exhibit a shelling of the boundary of the octahedron:
    \begin{itemize}
        \item[$F_1$]$=\{A,D,E\}$. In general, the restriction of $F_1$ is ${\rm R}(F_1)=\emptyset$.
        \item[$F_2$]$=\{A,D,F\}$. We have ${\rm R}(F_2)=\{F\}$.
        \item[$F_3$]$=\{C,D,E\}$. We have ${\rm R}(F_3)=\{C\}$.
        \item[$F_4$]$=\{A,B,E\}$. We have ${\rm R}(F_4)=\{B\}$. 
        \item[$F_5$]$=\{A,B,F\}$. We have ${\rm R}(F_5)=\{A,F\}$.
        \item[$F_6$]$=\{B,C,E\}$. We have ${\rm R}(F_6)=\{B,C\}$.
        \item[$F_7$]$=\{C,D,F\}$. We have ${\rm R}(F_7)=\{F,C\}$ .
        \item[$F_8$]$=\{B,C,F\}$. We have ${\rm R}(F_8)=\{B,C,F\}$.
    \end{itemize}
    By~\eqref{h-vec as restrictions}, this gives the $h$-vector $(1, 3, 3, 1)$.
    Applying \eqref{eq-h-to-f}, we retrieve the $f$-vector $(6,12,8)$. 
    Indeed, the boundary of octahedron contains $6$ vertices, $12$ edges, and $8$ faces (i.e., facets).
\end{example}

\subsection{The Stanley--Reisner ring}
\label{sub:SR}

Let $\Delta$ be a simplicial complex on the vertex set $V$.
Let $K$ be a field, and consider the polynomial ring $K[V] \coloneqq K[x_v : v \in V]$, where we regard the elements of $V$ as indeterminates.
Given a subset $U \subseteq V$, we will use the shorthand
\[
x_U \coloneqq \prod_{v \in U} x_v \qquad \text{and} \qquad K[U] \coloneqq K[x_v : v \in U ].
\]

Let $I_\Delta$ be the ideal of $K[V]$ generated by all monomials $x_U$ such that $U \not\in \Delta$.
Such a $U$ is called a \emph{nonface} of $\Delta$, and it is easy to see that $I_\Delta$ is actually generated by those nonfaces which are minimal (i.e., which contain no proper nonface).
The \emph{support} of a monomial $\mathbf m \in K[V]$ is the set $\{ v \in V : \text{$x_v$ divides $\mathbf m$}\}$.
A $K$-basis for $I_\Delta$ is given by the monomials whose support is not contained in $\Delta$.   

The quotient $K[\Delta] \coloneqq K[V]/I_\Delta$ is called the \emph{Stanley--Reisner} ring of $\Delta$.  
It is clear that $K[\Delta]$ has a $K$-basis consisting of the monomials whose support is a face of $\Delta$ (where we identify these monomials with their images in the quotient ring).

A shelling of $\Delta$ induces a \emph{Stanley decomposition} of the Stanley--Reisner ring:
\begin{equation}
    \label{Stanley decomp general}
    K[\Delta] = \bigoplus_{F} x_{{\rm R}(F)} K[F],
\end{equation}
where the direct sum ranges over all facets $F \in \Delta$, and where their restrictions ${\rm R}(F)$ are determined by the shelling.
Crucially, each monomial in $K[\Delta]$ lies in exactly one summand of~\eqref{Stanley decomp general}.
Since $I_\Delta$ is generated by homogeneous polynomials (in fact, by monomials), the quotient $K[\Delta]$ inherits from $K[V]$ the natural grading by degree.
Writing $K[\Delta]_{d}$ to denote the graded component consisting of homogeneous polynomials of degree $d$, we can restrict~\eqref{Stanley decomp general} to a decomposition of each component:
\begin{equation}
    \label{Stanley decomp graded}
    K[\Delta]_d = \bigoplus_k \bigoplus_{\substack{F:\\ \#{\rm R}(F) = k}} x_{{\rm R}(F)} K[F]_{d-k}.
\end{equation}
Clearly the inside summand is empty if $k > d$ or if $k$ exceeds the size of the largest restriction ${\rm R}(F)$; hence $k$ ranges from $0$ to the minimum of $d$ and this size.

\subsection{The order complex of $\Pi_{ij}$}
\label{sub:order complex}

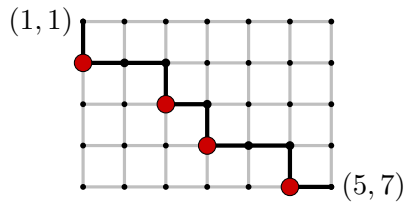
\begin{figure}[ht]
    \centering
    \tikzset{every picture/.style={scale=.55}}
\begin{tikzpicture}
\draw[gray!50, very thick] (0,0) grid (6,4);
\foreach \i in {0,...,6}
    \foreach \j in {0,...,4}
        \fill (\i,\j) circle (2pt);
\draw[ultra thick] (0,4) -- (0,3)--(2,3)--(2,2)--(3,2)--(3,1)--(5,1)--(5,0)--(6,0) ;
\fill[black]  
(1,3) circle (3pt)
(2,3) circle (3pt)
(2,2) circle (3pt) 
(3,2) circle (3pt) 
(3,1) circle (3pt) 
(4,1) circle (3pt) 
(5,1) circle (3pt) 
(5,0) circle (3pt);

\path [draw=black, fill=red!80!black]
(0,3) circle (6pt) 
(2,2) circle (6pt)
(3,1) circle (6pt)
(5,0) circle (6pt);           
\fill[black] 
(0,4) circle (2pt) node[left] {$(1,1)$} 
(6,0) circle (2pt) node[right] {$(5,7)$};
\end{tikzpicture}
    \caption{Example of a facet $F$ of $\Delta_{5,7}$, or equivalently, a maximal chain in $\Pi_{5,7}$.  
    The four red circles (the \emph{corners} of $F$) indicate the restriction ${\rm R}(F)$.  
    This facet thus contributes $t^4$ to the polynomial $h_{5,7}(t)$.}
    \label{fig:path example}
\end{figure}

We now apply the general theory above to the main problem of this paper, i.e., writing down a formula for the matrix $\SS(d,n)$.  
Recall the poset $\Pi_{ij}$ defined in~\eqref{def:Pi}.  
We define the \emph{order complex} $\Delta_{ij}$ on the vertex set $\Pi_{ij}$, such that the faces of $\Delta_{ij}$ are the chains in $\Pi_{ij}$.  
Then the facets of $\Delta_{ij}$ are the maximal chains $\{(1,1), \ldots, (i,j)\}$, all of which have size $i+j-1$; hence $\Delta_{ij}$ is pure.  
Note that each facet can be visualized as a lattice path from the northwest to the southeast corner of an $i \times j$ matrix, consisting of steps to the south and to the east; see Figure~\ref{fig:path example}.  

Since $\Pi_{ij}$ is a finite planar distributive lattice, it follows from~\cite{Bjorner}*{Thm.~7.1} that $\Pi_{ij}$ is shellable.
Specifically, if we label each step toward the east with ``1'' and each step toward the south with ``2,'' then we can identify each facet with a unique sequence of 1's (of which there are $j-1$) and 2's (of which there are $i-1$).
This induces a lexicographical order on the facets, which we take as our shelling.
With respect to this shelling, each restriction ${\rm R}(F)$ is given by the set of descents in the sequence of 1's and 2's corresponding to $F$ (where a descent is a 2 followed immediately by a 1).
Hence, a descent corresponds to a point $(x,y) \in F$ at which a south step ends and an east step begins:
\begin{equation}
    \label{corners}
    {\rm R}(F) = \{ (x,y) \in F : \text{both $(x-1,\: y) \in F$ and $(x,\: y+1) \in F$}\}.
\end{equation}
Visually, then, the elements in ${\rm R}(F)$ occur at the \scalebox{2}{$\llcorner$}-corners in the path representing $F$, as shown in Figure~\ref{fig:path example}.
For this reason we will sometimes refer to elements of ${\rm R}(F)$ as the \emph{corners} of $F$.
(The identification of the restrictions with sets of descents follows from the general theory of lexicographical shellings in~\cite{BjornerWachs}*{Thm.~5.8}.
See also~\cite{BCS}*{\S3} for a detailed treatment of this specific shelling.)  This realization of restrictions as sets of corners yields a simple formula for the $h$-polynomial of $\Delta_{ij}$:

\begin{lemma}
    \label{lemma:h-poly}
    Let $h_{ij}(t)$ be the $h$-polynomial of $\Delta_{ij}$.  Then we have 
    \[
    h_{ij}(t) = \sum_{\ell=0}^{\min\{i-1, \: j-1\}} \binom{i-1}{\ell}\binom{j-1}{\ell} t^\ell.
    \]
\end{lemma}

\begin{proof}
    Suppose that $F$ is a facet in $\Delta_{ij}$ such that $\#{\rm R}(F) = \ell$.  
    Then $F$ is a maximal path in $\Pi_{ij}$ with $\ell$ corners: 
    \[
    {\rm R}(F) = \{(x_1, y_1), \ldots, (x_\ell, y_\ell)\}, \qquad 2 \leq x_1 < \cdots < x_\ell \leq i \text{ and } 1 \leq y_1 < \cdots < y_\ell \leq j-1.
    \]
Hence there are $\binom{i-1}{\ell}$ choices for the $x$'s and $\binom{j-1}{\ell}$ choices for the $y$'s, each of which determines a unique set ${\rm R}(F)$ and therefore a unique facet $F$.  
The lemma follows from ~\eqref{h-vec as restrictions}.
\end{proof}

\begin{rem}
    \label{rem:hij}
    We can also write $h_{ij}(t) = \leftidx{_2}{F}{_1}\!\left[\genfrac{}{}{0pt}{1}{1-i, \: 1-j}{1}; t\right]$, although we will not need this fact.
\end{rem}

We now turn to our main problem: write down a formula for each entry in the matrix $\SS(d,n)$.
Let $\Delta_n \coloneqq \Delta_{nn}$.
By the general theory in Section~\ref{sub:SR}, a $K$-basis for $K[\Delta_n]$ is given by the monomials whose support is a chain in $\Pi_n$.
Restricting to the degree-$d$ graded component, we observe the bijection
\begin{align*}
\T(d,n) & \longrightarrow \text{$K$-basis of $K[\Delta_n]_d$},\\
T & \longmapsto \prod_{(i,j) \in \Pi_n} x_{ij}^{T_{ij}},
\end{align*}
obtained by assigning the exponent matrix $T$ to each monomial in $K[\Delta_n]$.
Under this correspondence, adding matrices corresponds to multiplying monomials.
It follows that, letting $\mathbf m$ range over the monomials,
\begin{equation}
    \label{adding arrays is mult monomials}
    \prod_{\mathclap{\mathbf m \in K[\Delta_n]_d}} \mathbf m = \prod_{(i,j) \in \Pi_n} x_{ij}^{\SS(d,n)_{ij}},
\end{equation}
i.e., $\SS(d,n)$ is the exponent matrix of the product of all monomials in $K[\Delta_n]_d$.
Hence our main problem is equivalent to writing down this product of monomials, which we carry out in the next section.

\begin{rem} 
The structure of our Stanley--Reisner ring $K[\Delta_n]$ is particularly well studied in several areas of mathematics.
First, it is isomorphic to the coordinate ring of the determinantal variety of $n \times n$ matrices with rank $\leq 1$, which was the motivation behind~\cite{BCS}.  
Lattice-path arguments like those in our paper have often been used in connection with determinantal varieties; see~\cites{Sturmfels,Herzog,Conca}, for example.
Second, if $K$ has characteristic $0$, then $K[\Delta_n]$ is also isomorphic to the invariant ring $K[x_1, \ldots, x_n,y_1, \ldots, y_n]^{K^\times}$, where $\alpha \in K^\times$ acts on the variables via $\alpha \cdot x_i = \alpha x_i$ and $\alpha \cdot y_i = \alpha^{-1}y_i$.  The explicit isomorphism is given by $x_{ij} \mapsto x_i y_j$.  
This is a special case of one of the central settings in classical invariant theory, namely, the general linear group ${\rm GL}(r,K)$ acting on an arbitrary number of vectors and covectors;
in general, the invariant ring is isomorphic to the coordinate ring of the determinantal variety of rank $\leq r$.
This is the substance of the celebrated first and second fundamental theorems of Weyl~\cite{Weyl}.
\end{rem}

\section{Main result, second version}
\label{sec:main result h-poly}

We present the second version of our main result:

\begin{theorem}
    \label{thm:main result h-polys}
    For all $1 \leq i,j \leq n$, the $(i,j)$ entry of $\SS(d,n)$ is given by
    \[
    \SS(d,n)_{ij} = \sum_{k=0}^{\min\{d,n\}-1} \binom{2n+d-k-2}{2n-1} \sum_{\ell = 0}^k \binom{i-1}{\ell}\binom{j-1}{\ell}\binom{n-i}{k-\ell}\binom{n-j}{k-\ell}.
    \]
\end{theorem}

Before giving the proof, we record two easy counting lemmas. 

\begin{lemma}
    \label{lemma:alpha}
    Let $(i,j) \in \Pi_n$, and let $F \ni (i,j)$ be a facet of $\Delta_n$.
    Then $\binom{2n+d-k-2}{2n-1}$ equals the exponent of $x_{ij}$ in the product of all monomials in
\begin{equation}
    \label{summand in lemma}
    x_{ij} x_{{\rm R}(F) \setminus \{(i,j)\}} \:K[F]_{d-k-1}.
\end{equation}
\end{lemma}

\begin{proof}
It suffices to show that
\[
\binom{2n+d-k-2}{2n-1} = \Big(\text{\# monomials in \eqref{summand in lemma}}\Big) \Big(\text{average exponent of $x_{ij}$ in each monomial}\Big).
\]
The number of monomials in~\eqref{summand in lemma} equals the number of monomials in $K[F]_{d-k-1}$, which is the number of weak compositions of the degree $d-k-1$ into $\#F$ many parts.
Thus, recalling that $\#F = 2n-1$ for any facet $F$ of $\Delta_n$, and using the elementary formula~\eqref{composition formula}, we have
\begin{equation}
    \label{number monomials}
    \text{\# monomials in~\eqref{summand in lemma}} = \binom{d-k-1 + (2n-1) - 1}{(2n-1)-1} = \binom{2n+d-k-3}{2n-2}.
\end{equation}
The average exponent of $x_{ij}$, taken over all the monomials in $K[F]_{d-k-1}$, equals the degree $d-k-1$ divided by the number $\#F$ of variables.
Adding 1 to this average to account for the factor of $x_{ij}$ present in~\eqref{summand in lemma}, we obtain 
\begin{equation}
    \label{avg exp}
    \text{average exponent of $x_{ij}$ in each monomial} = 1 + \frac{d-k-1}{2n-1} = \frac{2n + d - k - 2}{2n-1}.
\end{equation}
Multiplying the expressions in~\eqref{number monomials} and~\eqref{avg exp}, we obtain
\[
\binom{2n + d -k -3}{2n-2} \cdot \frac{2n+d-k-2}{2n-1} = \binom{2n+d-k-2}{2n-1},
\]
as desired.
\end{proof}

\begin{lemma}
    \label{lemma: coeff tk}
    Let $(i,j) \in \Pi_n$.
    Then 
    \[
    \sum_{\ell = 0}^k \binom{i-1}{\ell}\binom{j-1}{\ell}\binom{n-i}{k-\ell}\binom{n-j}{k-\ell}
    \]
    equals the number of facets $F \ni (i,j)$ of $\Delta_n$ such that $\#\Big({\rm R}(F) \setminus \{(i,j)\}\Big)= k$.
\end{lemma}

\begin{proof}
    Every facet $F \ni (i,j)$ of $\Delta_n$ is the union of two saturated chains
    \[
    F' : (1,1) \preceq \cdots \preceq (i,j) \qquad \text{and} \qquad F'' : (i,j) \preceq \cdots \preceq (n,n),
    \]
    which intersect only at $(i,j)$.
    Clearly $F'$ can be any facet of $\Delta_{ij}$, viewed as a subposet of $\Delta_n$.
    Likewise, $F''$ can be any facet of $\Delta_{n-i+1, n
    -j+1}$, viewed as a subposet of $\Delta_n$ after translating coordinates.
    Since $(i,j)$ is either the maximal or minimal element of these two subposets, it cannot occur as an element of ${\rm R}(F')$ or of ${\rm R}(F'')$.
    Therefore $\#({\rm R}(F) \setminus \{(i,j)\})  = {\rm R}(F') + {\rm R}(F'')$.
    By~\eqref{h-vec as restrictions}, we thus have
    \begin{align*} 
    h_{ij}(t) h_{n-i+1, n-j+1}(t) &= \left(\sum_{F'} t^{\#{\rm R}(F')} \right)\left(\sum_{F''} t^{\#{\rm R}(F'')}\right)\\
    &= \sum_{F', F''} t^{\#{\rm R}(F') + \#{\rm R}(F'')}\\
    &= \sum_{F \ni (i,j)} t^{\#({\rm R}(F) \setminus \{(i,j)\})},
    \end{align*}
    where the sums range over facets $F$, $F'$, and $F''$ of $\Delta_n$, $\Delta_{ij}$, and $\Delta_{n-i+1, n-j+1}$, respectively.
    Therefore the number of facets described in the lemma equals the coefficient of $t^k$ in the following product, which we expand via Lemma~\ref{lemma:h-poly}:
    \[
    h_{ij}(t) h_{n+1-i, n+1-j}(t) = \left(\sum_{\ell = 0}^{\min\{i-1,\:j-1\}} \binom{i-1}{\ell}\binom{j-1}{\ell} t^\ell \right)\left(\sum_{m = 0}^{\min\{n-i,\:n-j\}} \binom{n-i}{m}\binom{n-j}{m} t^m \right).
    \]
    The coefficient of $t^k$ in this expansion equals $\sum_{\ell + m = k} \binom{i-1}{\ell}\binom{j-1}{\ell} \binom{n-i}{m}\binom{n-j}{m}$.
    Upon substituting $k-\ell$ for $m$, the proof is complete.
\end{proof}

\begin{proof}[Proof of Theorem \ref{thm:main result h-polys}]

By~\eqref{Stanley decomp graded} and~\eqref{adding arrays is mult monomials}, we know that $\SS(d,n)_{ij}$ equals the exponent of $x_{ij}$ in the product of all monomials in the graded component
\begin{equation}
    \label{Stanley decomp in proof first}
    K[\Delta_n]_d = \bigoplus_{k=0}^{\min\{d, \:n-1\}} \bigoplus_{\substack{F:\\ \#{\rm R}(F) = k}} x_{{\rm R}(F)} K[F]_{d-k},
\end{equation}
where the inside sum ranges over the facets $F$ of $\Delta_n$.
But the only monomials contributing to this exponent are those divisible by $x_{ij}$.
Hence we may ignore all summands in~\eqref{Stanley decomp in proof first} such that $(i,j) \not\in F$.
If $(i,j) \in F$, then the subspace of $K[F]_{d-k}$ spanned by the monomials divisible by $x_{ij}$ is
\[
x_{ij} \:K[F]_{d-k-1}.
\]
Then since $(i,j)$ may or may not lie in ${\rm R}(F)$, the subspace of $x_{{\rm R}(F)} K[F]_{d-k}$ spanned by monomials divisible by $x_{ij}$ is
\[
x_{ij} x_{{\rm R}(F) \setminus \{(i,j)\}} K[F]_{d-k-1}.
\]
Combining this with~\eqref{Stanley decomp in proof first}, we conclude that $\SS(d,n)_{ij}$ equals the exponent of $x_{ij}$ in the product of all monomials in
\begin{equation}
    \label{Stanley decomp in proof second}
    \bigoplus_{k=0}^{\min\{d-1,\:n-1\}} \bigoplus_{\substack{F \ni (i,j):\\ \#({\rm R}(F) \setminus \{(i,j)\}) = k}}  x_{ij} x_{{\rm R}(F) \setminus \{(i,j)\}} K[F]_{d-k-1}.
\end{equation}
Applying Lemma~\ref{lemma:alpha} to~\eqref{Stanley decomp in proof second}, we see that the desired exponent of $x_{ij}$ equals
\begin{align*}
    \SS(d,n)_{ij} &= \sum_{k=0}^{\min\{n,d\}-1} \sum_{\substack{F \ni (i,j):\\ \#({\rm R}(F) \setminus \{(i,j)\}) = k}} \binom{2n+d-k-2}{2n-1}\\[2ex]
    &= \sum_{k=0}^{\min\{d,n\}-1} \binom{2n+d-k-2}{2n-1} 
    \cdot \#\Big\{F \ni (i,j) : \#({\rm R}(F) \setminus \{(i,j)\}) = k \Big\},
\end{align*}
and we have already computed the second factor in Lemma~\ref{lemma: coeff tk}.
\end{proof}

Analogously to Remark~\ref{rem:n1n2}, one can easily extend the second version of our main result to $n_1 \times n_2$ matrices, obtaining
\[
    \SS(d,n_1 \times n_2)= \sum_{k=0}^{\min\{d,n_1,n_2\}-1} \binom{n_1+n_2+d-k-2}{n_1+n_2-1} \sum_{\ell = 0}^k \binom{i-1}{\ell}\binom{j-1}{\ell}\binom{n_1-i}{k-\ell}\binom{n_2-j}{k-\ell}.
\]

\begin{rem}
    In the formula in Theorem~\ref{thm:main result h-polys}, it is certainly possible to express the sum over $\ell$ (i.e., the result in Lemma~\ref{lemma: coeff tk}) in terms of a $\leftidx{_4}{F}{_3}$ series, using the identity~\eqref{sum of 2F1s} with the fact that $h_{ij}(t)$ is a $\leftidx{_2}{F}{_1}$ series (Remark~\ref{rem:hij}).
    Unlike in our first formula, however, we have yet to find an expression that closes the sum over $k$.
    \end{rem}

\section{Application to optimal transport}
\label{sec:EMD}

We conclude by highlighting one application of our main result to optimal transport theory.
Since its 18th-century origin in a treatise of Gaspard Monge~\cite{monge}, and its modern re-popularization by Rubner et.\ al.\ \cite{rubner} in computer vision, the classical transportation problem --- more precisely, its solution, known as the earth mover's distance (EMD) --- has become central to applications in physics \cite{komiske}, cosmology \cite{frisch}, political science \cite{lupu}, epidemiology \cite{melnyk}, and many others; see~\cite{Villani} for a comprehensive treatment.
Below, we consider the discrete setting of the transportation problem.

Let $\CC(d,n)$ denote the set of weak integer compositions of $d$ into $n$ parts.
In particular, let $\la = (\la_1, \ldots, \la_n) \in \CC(d,n)$ be the ``supply vector,'' where $\la_i$ gives the number of units at supply location $i$; likewise, let $\mu = (\mu_1, \ldots, \mu_n)$ be the ``demand vector,'' where $\mu_j$ gives the number of units required at the demand location $j$.
For $1 \leq i,j \leq n$, let $C_{ij} \geq 0$ denote the cost of transporting one unit from supply location $i$ to demand location $j$.
The goal, of course, is to transport all units from the supply locations to the demand locations, so as to minimize the total cost.  
In particular, we want to find an $n \times n$ matrix $T$ in order to solve the following linear programming problem:
\begin{align*}
    \text{Minimize} \quad & \sum_{\mathclap{i,j=1}}^n C_{ij} T_{ij},\\
    \text{subject to} \quad & T_{ij} \geq 0 & \text{for all } 1 \leq i,j \leq n,\\
    \text{and} \quad & \sum_{j=1}^n T_{ij} = \la_i & \text{for each }1\leq i \leq n,\\
    \text{and} \quad & \sum_{i=1}^n T_{ij} = \mu_j & \text{for each } 1 \leq j \leq n.
\end{align*}
Thus one can regard $T$ as a contingency table with $\la$ and $\mu$ as its margins.
The solution to the transportation problem above is nicest when the cost matrix $C$ has what is known as the \emph{Monge property}:
\begin{equation}
    \label{Monge property}
    C_{ij} + C_{IJ} \leq C_{Ij} + C_{iJ}, \quad \text{for all $i<I$ and $j < J$.}
\end{equation}
In this case, as shown by Hoffman~\cite{hoffman}, the problem can be solved via a certain $O(2n)$-time greedy algorithm called the \emph{northwest corner rule}.  
The northwest corner rule solves the transportation problem by constructing an optimal transport matrix $T$ as follows: 
beginning in the upper-left, the algorithm sets $T_{11} = \min\{\la_1,\mu_1\}$, and then modifies both $\la_1$ and $\mu_1$ by subtracting $T_{11}$ (so at least one of them becomes 0).  
This allows us to fill the remainder of either the first row or the first column with 0's;
then we proceed either south or east, and repeat the process until we have filled the entire matrix $T$, which is the solution to the transportation problem. 

It is immediate from this construction that $T$ has nonnegative integer entries summing to $d$, and has row and column sums given by $\la$ and $\mu$, respectively.  
Moreover, the very nature of the northwest corner rule forces the nonzero entries of $T$ to lie in a single path traveling either south or east at each step: 
in other words, $T$ is width-one and hence $T \in \T(d,n)$, as defined in~\eqref{T(d,n) definition}. 
Therefore, writing $T_{\la\mu}$ for the output of the northwest corner rule applied to the pair $(\la,\mu$), we have a bijection
\begin{align}
\label{bijection histogram pairs CTs}
\begin{split}
 \CC(d, n) \times \CC(d,n) &\xrightarrow{\phantom{\text{NW corner rule}}}\T(d,n),\\
 (\la,\mu) &\xmapsto{\text{NW corner rule}} T_{\la\mu},
 \end{split}
\end{align}
which is essentially our bijection from Section~\ref{sec:RSK} between biwords and width-one matrices.

In the literature, the minimum value of the objective function $\sum_{ij} C_{ij} T_{ij}$ above is often called the ``earth mover's distance'' (EMD) between $\la$ and $\mu$.  
In \cite{BW}, Bourn and Willenbring derived a recursive formula to compute the mean value of the EMD, in the context of one-dimensional histograms with bins located at consecutive points $1, \ldots, n$ on the number line. 
In this case the cost of transport between two bins is the $\ell_1$-distance between them.  
Each author of the present paper has worked independently to generalize the results of~\cite{BW} in different directions \cites{EricksonAStat, KretschmannMasters, KretschmannPhD}, which led us to seek an explicit, non-recursive formula for the mean EMD between discrete histogram pairs.
(We point out that the paper~\cite{FV} took an analytical approach to close the recursion in~\cite{BW} in the setting of probability distributions rather than discrete histograms.)
It turns out that this problem is easily solved by applying our main result in this paper:

\begin{prop}
    \label{prop:EMD=tr(CS)}
    Let the cost matrix $C$ be an $n \times n$ matrix with the Monge property~\eqref{Monge property}.
    Then the mean EMD on $\CC(d,n) \times \CC(d,n)$ equals
    \[
    \frac{1}{\binom{d+n-1}{d}^2} \cdot \operatorname{tr}\!\Big(C^\mathsf{T}\! \cdot \SS(d,n)\Big).
    \]
    \end{prop}

\begin{proof}

To obtain the mean, we multiply the reciprocal of $\#\CC(d,n)^2 = \binom{d+n-1}{d}^2$ by the following sum, which we rewrite using the linearity of the trace:
\begin{align*}
    \sum_{\mathclap{\substack{(\la,\mu) \\ \in \CC(d,n) \times \CC(d,n)}}} {\rm EMD}(\la,\mu) &= \sum_{(\la,\mu)} \sum_{i,j=1}^n C_{ij}(T_{\la\mu})_{ij}\\
    &= \sum_{(\la,\mu)} \operatorname{tr}(C^\mathsf{T}  T_{\la\mu})\\
    &= \operatorname{tr}\left(C^\mathsf{T} \! \cdot \sum_{(\la,\mu)} T_{\la\mu}\right)\\
    &= \operatorname{tr}\left(C^\mathsf{T} \! \cdot \sum_{\mathclap{T \in \T(d,n)}} T\right)\\
    &= \operatorname{tr}\Big(C^\mathsf{T} \! \cdot \SS(d,n)\Big),
\end{align*}
where the last two equalities follow from \eqref{bijection histogram pairs CTs} and \eqref{Bold T definition}, respectively.
\end{proof}

Hence our explicit formulas for $\SS(d,n)$ allow one to compute the mean EMD by taking the trace of a single matrix product.  
This is not only more direct, but also much more flexible than the recursive methods mentioned above: one can now immediately compute the mean EMD with respect to any cost matrix $C$ with the Monge property.

\begin{appendices}
    \section{Some examples of $\SS(d,n)$}
    \label{sec:appValues}

In this appendix, we display the matrices $\SS(d,5)$ for the first few values of $d$:
    \begin{align*}
\SS(1,5)&=
\left(\begin{array}{rrrrr}
1 & 1 & 1 & 1 & 1 \\
1 & 1 & 1 & 1 & 1 \\
1 & 1 & 1 & 1 & 1 \\
1 & 1 & 1 & 1 & 1 \\
1 & 1 & 1 & 1 & 1
\end{array}\right)
\\
\SS(2,5)&=
\left(\begin{array}{rrrrr}
26 & 22 & 18 & 14 & 10 \\
22 & 20 & 18 & 16 & 14 \\
18 & 18 & 18 & 18 & 18 \\
14 & 16 & 18 & 20 & 22 \\
10 & 14 & 18 & 22 & 26
\end{array}\right)
\\
\SS(3,5)&=
\left(\begin{array}{rrrrr}
251 & 193 & 141 & 95 & 55 \\
193 & 173 & 150 & 124 & 95 \\
141 & 150 & 153 & 150 & 141 \\
95 & 124 & 150 & 173 & 193 \\
55 & 95 & 141 & 193 & 251
\end{array}\right)
\\
\SS(4,5)&=
\left(\begin{array}{rrrrr}
1476 & 1064 & 720 & 440 & 220 \\
1064 & 960 & 816 & 640 & 440 \\
720 & 816 & 848 & 816 & 720 \\
440 & 640 & 816 & 960 & 1064 \\
220 & 440 & 720 & 1064 & 1476
\end{array}\right)
\\
\SS(5,5)&=
\left(\begin{array}{rrrrr}
6376 & 4385 & 2805 & 1595 & 715 \\
4385 & 4006 & 3360 & 2530 & 1595 \\
2805 & 3360 & 3546 & 3360 & 2805 \\
1595 & 2530 & 3360 & 4006 & 4385 \\
715 & 1595 & 2805 & 4385 & 6376
\end{array}\right)
\\
\SS(6,5)&=
\left(\begin{array}{rrrrr}
22252 & 14762 & 9042 & 4862 & 2002 \\
14762 & 13672 & 11352 & 8272 & 4862 \\
9042 & 11352 & 12132 & 11352 & 9042 \\
4862 & 8272 & 11352 & 13672 & 14762 \\
2002 & 4862 & 9042 & 14762 & 22252
\end{array}\right)
\\
\SS(7,5)&=
\left(\begin{array}{rrrrr}
66352 & 42779 & 25311 & 13013 & 5005 \\
42779 & 40150 & 33066 & 23452 & 13013 \\
25311 & 33066 & 35706 & 33066 & 25311 \\
13013 & 23452 & 33066 & 40150 & 42779 \\
5005 & 13013 & 25311 & 42779 & 66352
\end{array}\right)
\\
\SS(8,5)&=
\left(\begin{array}{rrrrr}
175252 & 110396 & 63492 & 31460 & 11440 \\
110396 & 104896 & 85800 & 59488 & 31460 \\
63492 & 85800 & 93456 & 85800 & 63492 \\
31460 & 59488 & 85800 & 104896 & 110396 \\
11440 & 31460 & 63492 & 110396 & 175252
\end{array}\right)
\\
\end{align*}
\end{appendices}

 \bibliographystyle{plain}
 \bibliography{references}

\begin{bibdiv}
\begin{biblist}

\bib{BCS}{article}{
      author={Billera, Louis~J.},
      author={Cushman, Richard~H.},
      author={Sanders, Jan~A.},
       title={The {S}tanley decomposition of the harmonic oscillator},
        date={1988},
        ISSN={0019-3577},
     journal={Nederl. Akad. Wetensch. Indag. Math.},
      volume={50},
      number={4},
       pages={375\ndash 393},
      review={\MR{976522}},
}

\bib{Bjorner}{article}{
      author={Bj\"{o}rner, Anders},
       title={Shellable and {C}ohen-{M}acaulay partially ordered sets},
        date={1980},
        ISSN={0002-9947},
     journal={Trans. Amer. Math. Soc.},
      volume={260},
      number={1},
       pages={159\ndash 183},
         url={https://doi.org/10.2307/1999881},
      review={\MR{570784}},
}

\bib{BjornerWachs}{article}{
      author={Bj\"{o}rner, Anders},
      author={Wachs, Michelle~L.},
       title={Shellable nonpure complexes and posets. {I}},
        date={1996},
        ISSN={0002-9947},
     journal={Trans. Amer. Math. Soc.},
      volume={348},
      number={4},
       pages={1299\ndash 1327},
         url={https://doi.org/10.1090/S0002-9947-96-01534-6},
      review={\MR{1333388}},
}

\bib{BW}{article}{
      author={Bourn, Rebecca},
      author={Willenbring, Jeb~F.},
       title={Expected value of the one-dimensional earth mover's distance},
        date={2020},
     journal={Algebr. Stat.},
      volume={11},
      number={1},
       pages={53\ndash 78},
}

\bib{Herzog}{article}{
      author={Bruns, Winfried},
      author={Herzog, J\"{u}rgen},
       title={On the computation of {$a$}-invariants},
        date={1992},
        ISSN={0025-2611},
     journal={Manuscripta Math.},
      volume={77},
      number={2-3},
       pages={201\ndash 213},
         url={https://doi.org/10.1007/BF02567054},
      review={\MR{1188581}},
}

\bib{Conca}{article}{
      author={Conca, Aldo},
       title={Ladder determinantal rings},
        date={1995},
        ISSN={0022-4049},
     journal={J. Pure Appl. Algebra},
      volume={98},
      number={2},
       pages={119\ndash 134},
         url={https://doi.org/10.1016/0022-4049(94)00039-L},
      review={\MR{1319965}},
}

\bib{Erdelyi}{book}{
      author={Erd\'{e}lyi, Arthur},
      author={Magnus, Wilhelm},
      author={Oberhettinger, Fritz},
      author={Tricomi, Francesco~G.},
       title={Higher transcendental functions. {V}ol. {I}},
   publisher={Robert E. Krieger Publishing Co., Inc., Melbourne, Fla.},
        date={1981},
        ISBN={0-89874-069-X},
      review={\MR{698779}},
}

\bib{EricksonAStat}{article}{
      author={Erickson, William~Q.},
       title={A generalization for the expected value of the earth mover's
  distance},
        date={2021},
        ISSN={2693-2997},
     journal={Algebr. Stat.},
      volume={12},
      number={2},
       pages={139\ndash 166},
         url={https://doi.org/10.2140/astat.2021.12.139},
      review={\MR{4350874}},
}

\bib{frisch}{article}{
      author={Frisch, Uriel},
      author={Matarrese, Sabino},
      author={Mohayaee, Roya},
      author={Sobolevski, Andrei},
       title={A reconstruction of the initial conditions of the universe by
  optimal mass transportation},
        date={2002},
     journal={Nature},
      volume={417},
       pages={260–262},
        note={https://doi.org/10.1038/417260a},
}

\bib{FV}{article}{
      author={Frohmader, Andrew},
      author={Volkmer, Hans},
       title={1-{W}asserstein distance on the standard simplex},
        date={2021},
        ISSN={2693-2997},
     journal={Algebr. Stat.},
      volume={12},
      number={1},
       pages={43\ndash 56},
         url={https://doi.org/10.2140/astat.2021.12.43},
      review={\MR{4251253}},
}

\bib{hoffman}{incollection}{
      author={Hoffman, A.J.},
       title={On simple linear programming problems},
        date={1963},
   booktitle={{Convexity: Proceedings of the Seventh Symposium in Pure
  Mathematics of the AMS}},
      editor={Klee, V.},
   publisher={American Mathematical Society},
     address={Providence, RI},
       pages={317\ndash 327},
}

\bib{24Hours}{book}{
      author={Iyengar, Srikanth~B.},
      author={Leuschke, Graham~J.},
      author={Leykin, Anton},
      author={Miller, Claudia},
      author={Miller, Ezra},
      author={Singh, Anurag~K.},
      author={Walther, Uli},
       title={Twenty-four hours of local cohomology},
      series={Graduate Studies in Mathematics},
   publisher={American Mathematical Society, Providence, RI},
        date={2007},
      volume={87},
        ISBN={978-0-8218-4126-6},
         url={https://doi.org/10.1090/gsm/087},
      review={\MR{2355715}},
}

\bib{Knuth}{article}{
      author={Knuth, Donald~E.},
       title={Permutations, matrices, and generalized {Y}oung tableaux},
        date={1970},
        ISSN={0030-8730},
     journal={Pacific J. Math.},
      volume={34},
       pages={709\ndash 727},
         url={http://projecteuclid.org/euclid.pjm/1102971948},
      review={\MR{272654}},
}

\bib{komiske}{article}{
      author={Komiske, Patrick~T.},
      author={Metodiev, Eric~M.},
      author={Thaler, Jesse},
       title={Metric space of collider events},
        date={2019},
     journal={Phys. Rev. Lett.},
      volume={123},
       pages={041801},
         url={https://link.aps.org/doi/10.1103/PhysRevLett.123.041801},
}

\bib{KretschmannPhD}{thesis}{
      author={Kretschmann, Jan},
       title={Title {TBD}},
        type={Ph.D. Thesis},
 institution={University of Wisconsin--Milwaukee},
        note={Expected 2024},
}

\bib{KretschmannMasters}{thesis}{
      author={Kretschmann, Jan},
       title={Earth mover's distance between grade distribution data with fixed
  mean},
        type={Master's Thesis},
 institution={University of Wisconsin--Milwaukee},
        date={2020},
         url={https://dc.uwm.edu/etd/2542},
}

\bib{lupu}{article}{
      author={Lupu, Noam},
      author={Selios, Luc\'ia},
      author={Warner, Zach},
       title={A new measure of congruence: the earth mover's distance},
        date={2017},
     journal={Political Analysis},
      volume={25},
       pages={95–113},
}

\bib{melnyk}{article}{
      author={Melnyk, Andrew},
      author={Knyazev, Sergey},
      author={Vannberg, Fredrik},
      author={Bunimovich, Leonid},
      author={Skums, Pavel},
      author={Zelikovsky, Alex},
       title={Using earth mover’s distance for viral outbreak
  investigations},
        date={2020},
     journal={BMC Genomics},
      volume={5},
       pages={582\ndash 590},
        note={https://doi.org/10.1186/s12864-020-06982-4},
}

\bib{monge}{incollection}{
      author={Monge, Gaspard},
       title={M\'emoire sur la th\'eorie des d\'eblais et des remblais},
        date={1781},
   booktitle={Histoire de l'acad\'emie royale des sciences de paris},
       pages={666\ndash 704},
}

\bib{PWZ}{book}{
      author={Petkov\v{s}ek, Marko},
      author={Wilf, Herbert~S.},
      author={Zeilberger, Doron},
       title={{$A=B$}},
   publisher={A.K. Peters, Ltd., Wellesley, MA},
        date={1996},
        ISBN={1-56881-063-6},
        note={With a foreword by Donald E. Knuth},
      review={\MR{1379802}},
}

\bib{rubner}{article}{
      author={Rubner, Yossi},
      author={Tomasi, Carlo},
      author={Guibas, Leonidas~J.},
       title={The earth mover's distance as a metric for image retrieval},
        date={2000},
     journal={Int. J. Comput. Vis.},
      volume={40},
       pages={99\ndash 121},
         url={https://doi.org/10.1023/A:1026543900054},
}

\bib{StanleyAC}{book}{
      author={Stanley, Richard~P.},
       title={Algebraic combinatorics: walks, trees, tableaux, and more},
      series={Undergraduate Texts in Mathematics},
   publisher={Springer, Cham},
        date={2018},
        ISBN={978-3-319-77172-4; 978-3-319-77173-1},
         url={https://doi.org/10.1007/978-3-319-77173-1},
        note={2nd ed.},
      review={\MR{3823204}},
}

\bib{Sturmfels}{article}{
      author={Sturmfels, Bernd},
       title={Gr\"{o}bner bases and {S}tanley decompositions of determinantal
  rings},
        date={1990},
        ISSN={0025-5874},
     journal={Math. Z.},
      volume={205},
      number={1},
       pages={137\ndash 144},
         url={https://doi.org/10.1007/BF02571229},
      review={\MR{1069489}},
}

\bib{Villani}{book}{
      author={Villani, C\'{e}dric},
       title={Optimal transport},
      series={Grundlehren der mathematischen Wissenschaften [Fundamental
  Principles of Mathematical Sciences]},
   publisher={Springer-Verlag, Berlin},
        date={2009},
      volume={338},
        ISBN={978-3-540-71049-3},
         url={https://doi.org/10.1007/978-3-540-71050-9},
        note={Old and new},
      review={\MR{2459454}},
}

\bib{Weyl}{book}{
      author={Weyl, Hermann},
       title={The classical groups: their invariants and representations},
   publisher={Princeton University Press, Princeton, NJ},
        date={1946},
      review={\MR{1488158}},
}

\end{biblist}
\end{bibdiv}

\end{document}